\documentclass[11pt]{article}
\usepackage{amsmath,amsfonts,amssymb}
\newtheorem{theorem}{Theorem}
\newtheorem{lemma}{Lemma}
\newtheorem{remark}{Remark}

\newcommand{\Z}{\mathbb{Z}}
\newcommand{\R}{\mathbb{R}}

\newcommand{\N}{{\mathbb N}}

\newcommand{\Fc}{\mathcal{F}}

\newcommand{\E}{\mathop{\mathsf E}}

\newcommand{\Pp}{\mathsf{P}}

\newcommand{\eps}{\varepsilon}

\newcommand{\ONE}{\mathbf{1}}

\newcommand{\bpf}[1][Proof]{{\noindent {\sc #1: }}}
\newcommand{\epf}{{{\hspace{4 ex} $\Box$ \smallskip}}}

\newcommand{\diam}{{\rm diam}}
\newcommand{\wb}{\overline}

\title{Localization and Perron--Frobenius
Theory for Directed Polymers}
\author{Yuri Bakhtin\thanks{
School of Mathematics, Georgia Institute of Technology, Atlanta, GA, 30067-0160, USA,
  email:{ bakhtin@math.gatech.edu}, phone: +1(404)894-9235, fax: +1(404)894-4409 (corresponding author)
} \and Konstantin Khanin\thanks{Department of Mathematics, University of Toronto, 40 St. George Street
Toronto, Ontario, M5S 2E4, Canada}}
\date{}
\begin{document}
\maketitle
\begin{abstract}
We consider directed polymers in a random potential given 
by a deterministic profile with a strong maximum at the origin 
taken with random sign at each integer time.

We study two main objects based on paths in this random potential.
First, we use the random potential and averaging over paths to define
a parabolic model via a random Feynman--Kac evolution operator. We
show that  for the resulting cocycle, there is a unique positive
cocycle eigenfunction serving as a forward and pullback attractor.
Secondly, we use the potential to define a Gibbs specification on
paths for any bounded time interval in the usual way and study the
thermodynamic limit and existence and uniqueness of an infinite
volume Gibbs measure.
Both main results claim that the local structure of interaction leads
to a unique macroscopic object for almost every realization of the
random potential.
\end{abstract}

\section{Introduction}
In this note, we consider two problems related to directed polymers with a localization property.
Directed $d+1$-dimensional polymers are modeled by random walks in a random potential 
$\phi:\Z^d\times\Z\times \Omega\to\R$
defined for all times $n\in\Z$ on the lattice $\Z^d$ via
\begin{equation}
\phi_n(x)=\phi_n(x,\omega)=V(x)B_n(\omega).
\label{eq:definition_of_potential}
\end{equation}
Here $V:\Z^d\to\R$ is a deterministic bounded function, and $B_n(\omega)$ is a sequence of
i.i.d. random variables defined on a probability space $(\Omega,\Fc,\Pp)$ and 
taking values $\pm1$ with probabilities $1/2$. It is convenient
to suppose that $\Omega=\{-1,1\}^\Z$, $\Fc$ is the cylindric $\sigma$-algebra, and $\Pp$ is the
$1/2$-Bernoulli product measure on $\Omega$. The random variables $B_n$ coincide with coordinate
maps: $B_n(\omega)=\omega_n$ for $\omega\in\Omega, n\in\Z$. 

For any integers $n_1\le n_2$ the set of {\it paths} that we consider consists of all possible trajectories
of the so called ``lazy random walk''. Namely, we deal only with paths
 $\gamma:[n_1,n_2]\cap\Z\to \Z^d$ such that
for any $n$, $\gamma_{n+1}$ equals either~$\gamma_n$ or one of the~$2d$ neighbors of $\gamma_n$ in $\Z^d$, i.e., $|\gamma_{n+1}-\gamma_n|\le 1$.

We 
study two main objects based on paths in the random potential $\phi$. First, we
use the random potential~$\phi$ and averaging over paths to define a parabolic model via a random Feynman--Kac evolution operator. We show that 
for the resulting cocycle, there is a unique positive cocycle eigenfunction serving as a forward and pullback attractor. 
Secondly, we use~$\phi$ to define a Gibbs specification on paths for any bounded time interval in the
usual way and study the thermodynamic limit and a.s.-existence and uniqueness of an infinite volume Gibbs measure.  
Both main results claim that the local structure of interaction leads to a unique macroscopic object for almost every
realization of the random potential.

The 
analogues of our results for a setting with a compact spin space replacing $\Z^d$ are based entirely on elementary
contraction arguments. The non-compactness of our spin space $\Z^d$ is the crucial feature that makes the problem
more interesting and difficult. We have to combine the contraction arguments with those based on localization. It is the
localization that makes the evolution of the system essentially compact.
The localization arguments we use are similar to those used by Sinai in his paper~\cite{SinaiMR1317991}
on Nechaev's model. However, the mechanism of localization in Nechaev's model is different from ours.
The result on existence and uniqueness of a Gibbs measure is also related to
the characterization of stationary Markov chains with finitely many states as Gibbs distributions with zero-range potential, see e.g.~\cite[Chapter 3]{GeorgiiMR956646}.

Localization for random walks in a random potential has been studied intensively in the past 10 years, see 
\cite{Carmona-Hu:MR1939654}, \cite{Comets-et-al:MR1996276}, \cite{Comets-et-al:MR2073332}, \cite{Alexander-Sidoravicius:MR2244428},
\cite{Giacomin:MR2380992},
and references therein. In our setting, the localization properties are based on a specific structure of the spatial
potential $V$. We assume that it has a large maximum (or minimum) located (without loss of generality) at the origin and has the
following form: $V(x)=V_0(x)+\Lambda \delta_0(x)$. Here $V_0$ is an arbitrary bounded potential and the delta-potential
$\delta_0$ takes value 1 at the origin and 0 everywhere else. All the results in this paper are proven under the assumption that
$\Lambda$ is large, the exact condition being expressed in terms of $\max|V_0(x)|$ and the dimension $d$. 
A special case when $V_0(x)\equiv 0$ was considered in \cite{Alexander-Sidoravicius:MR2244428},
\cite{Giacomin:MR2380992}.
It turn out that in dimension one and two localization occurs
for an arbitrary $\Lambda>0$.

In the theory of directed polymers it is often assumed that the dependence of the random potential on the space variable $x$
is non-deterministic. Of course, we can also think that $V_0$ is a realization of a bounded stationary potential $V_0^\omega$. However,
since our results hold for all such realizations we fix a particular realization and assume that $V$ is deterministic.

In the physics literature, mostly the model with independent noise values~$B_{x,n}$ for all space-time points $(x,n)$ has been considered. 
In the model that we suggest it is essential that one realization of the noise serves all lattice points $x$ at once, thus
modeling a spatially disordered structure represented by $V$ and embedded into an external potential $B_n$ fluctuating in time. 
The conditions we impose on the random potential~$\phi$ can be significantly relaxed though, and the results we prove hold true for a wider class of potentials with localization properties. 

For example, one can consider the following random potential:
\[
 \phi_n(x,\omega)=\sum_{i=1}^{m} V_i(x)B_n^{(i)}(\omega).
\]
Here $(V_i)_{1\le i\le m}$ is a collection of bounded potentials and $(B^{(i)})_{1\le i\le m}$ is a family of $m$ independent copies of the ``random sign'' process $B_n(\omega)$. To guarantee the localization it is sufficient to assume that one of the potentials~$V_i$ has a large
maximum or minimum at the origin. 

It is also possible to study the continuous time case where $B_n(\omega)$ is replaced by the white noise
$\dot W(t)$. However, this case is technically more complicated and it will be considered in a forthcoming publication. 

The paper is organized as follows. In this section, we describe the setting in detail and state our two main results. 
The crucial localization lemma is proven in Section~\ref{sec:main_loc_lemma}, and it is used to prove the main results in
Sections~\ref{sec:cocycle} (cocycle eigenfunctions) and~\ref{sec:Gibbs} (infinite volume Gibbs distributions) respectively.
Section~\ref{sec:auxi_lemmas} is devoted to auxiliary technical lemmas. 

\subsection{The parabolic model}

The Feynman--Kac operator associated with $\phi$ is defined for two integers $n_1<n_2$ and all
bounded functions $u$ as
\begin{align*}
T^{n_1,n_2}u(x)=(T^{n_1,n_2}(\omega) u)(x)=\E e^{\Phi_{n_1,n_2}(\gamma)} u(\gamma(n_1)),
\end{align*}
where 
\begin{equation*}
\Phi_{n_1,n_2}(\gamma)=\Phi_{n_1,n_2}(\gamma,\omega)=\sum_{n=n_1+1}^{n_2} \phi_n(\gamma(n),\omega)
=\sum_{n=n_1+1}^{n_2} V(\gamma(n))\omega_n,
\end{equation*} 
and the expectation $\E$ is taken with respect to the uniform distribution on the space of
directed polymer realizations $\Gamma_{n_1,n_2}(x)$ which is 
the set of admissible paths $\gamma:[n_1,n_2]\to\Z^d$ with $\gamma(n_2)=x$, i.e.
\begin{equation}
T^{n_1,n_2}u(x)=
\frac{1}{(2d+1)^{n_2-n_1}}\sum_{\gamma\in\Gamma_{n_1,n_2}(x)} e^{\Phi_{n_1,n_2}(\gamma)} u(\gamma(n_1)). 
\label{eq:T}
\end{equation}
We recall that throughout the paper, by {\it admissible paths} or, simply, {\it paths} we mean lazy random walk trajectories (see the
Introduction), and by 
$[n_1,n_2]$ we always mean a discrete set $\{n_1,\ldots,n_2\}$.


It is easy to see that $T^{n_1,n_2}$ can be considered as a discrete time analogue of the Feynman-Kac kernel for
the linear heat equation.
The family of operators $T^{n_1,n_2}$ defines a random dynamical system (cocycle) on $L^\infty=L^\infty(\Z^d)$
(we shall denote the $L^\infty$ norm by $\|\cdot\|$):
for every $\omega\in\Omega$
\begin{equation*}
T^{n_2,n_3}(\omega)T^{n_1,n_2}(\omega)=T^{n_1,n_3}(\omega),
\end{equation*} 
or, equivalently,
\[
T^{n_2}(\theta^{n_1}\omega)T^{n_1}(\omega)=T^{n_1+n_2}(\omega), 
\]
where $T^n(\omega)=T^{0,n}(\omega)$ and the time shift $\theta^k:\Omega\to\Omega$ is defined via
\begin{equation*}
\theta^k \omega(n)=\omega(k+n),\quad k,n\in\Z.
\end{equation*}

Our goal is to prove a Perron--Frobenius theorem for the cocycle $T$ defined above. The classical notion
of an eigenfunction has to be modified for the cocycle setting.

We say that an $L^\infty$-valued random variable $u^\omega$ is an eigenfunction for the cocycle $T$
if with probability 1 and for all $n\in\N$, $T^n(\omega)u^\omega$ is equal to~$u^{\theta^n\omega}$ up to
a random scalar factor. It is sufficient to impose this requirement for $n=1$. Obviously, $u^\omega$ is defined up to a multiplication by an arbitrary $c(\omega)\ne0$. Thus, it
is convenient to introduce normalized eigenfunctions assuming $\|u^\omega\|=1$. Then
\[
 T^n(\omega)u^{\omega}=\kappa_n(\omega)u^{\theta^n\omega},\quad n\in\N.
\]
for some process $(\kappa_n(\omega))_{n\in\N}$. When working with normalized functions, it is useful to consider
a normalized version of the cocycle $T^n$ defined by
\[
 \wb T^n(\omega) \phi =\frac{T^n(\omega) \phi}{\|T^n(\omega) \phi\|}.
\]

For a positive $\lambda$, an eigenfunction $u^\omega(x)$ is called $\lambda$-localized if there is a random variable $c=c(\omega)>0$ such that with probability 1,
for all $x\in\Z^d$,
\[
 |u^\omega(x)|\le c(\omega)e^{-\lambda|x|},
\]
where $|x|=\max\{|x_1|,\ldots,|x_d|\}$.


%

Throughout the paper we shall require the following conditions on the potential function $V$:
there are constants $M_0,M_1$ and $\lambda_0$
such that 
\begin{align}
\label{eq:conditions_on_V_1}
&V(0)=M_0,
\\
&|V(x)|\le M_1,\quad x\ne 0,
\label{eq:conditions_on_V_2}\\
&\lambda_0{=}\frac12 (M_0-3M_1)-\ln(2d+1)>0,
\label{eq:conditions_on_V_3}
\end{align}
For our first main theorem we shall need one more condition:
\begin{equation}
\lambda_1 \stackrel{def}{=} 2M_1+\ln(2d+1)<\lambda_0
\label{eq:conditions_on_V_4}
\end{equation}
\begin{theorem}\label{thm:perron-frobenius-cocycle} 
Suppose conditions \eqref{eq:conditions_on_V_1}---\eqref{eq:conditions_on_V_4} on $V$ are satisfied.
Then there is a unique normalized,  positive  eigenfunction~$u^\omega$ for the cocycle $T$ satisfying
\begin{equation}
\label{eq:uniqueness_condition}
\limsup_{n\to-\infty}u^{\theta^n\omega}(0)>0.
\end{equation}
This eigenfunction is $\lambda$-localized for any $\lambda\in(0,\lambda_0)$.

Moreover, almost surely,
for any nonnegative function $v\in L^\infty$ not identically equal to zero,
\begin{equation*}
\lim_{n\to\infty}\left\|{u^{\theta^n\omega}}-\wb T^{n}v\right\|= 0\quad\text{\rm (forward attraction),}
\end{equation*}
and
\begin{equation*}
\lim_{n\to\infty}\wb T^{-n,0}v= u^\omega\quad\text{\rm (pullback attraction)}.
\end{equation*}
\end{theorem}
A proof of this result is given in Section~\ref{sec:cocycle}.

\begin{remark} It follows straightforwardly from the ergodicity of the shift operator $\theta$ that there is a non-random Lyapunov
exponent $\lambda_L$ so that with probability 1,
\[ 
 \lim_{n\to\infty} \frac{\ln \kappa_n(\omega)}{n}=\lambda_L.
\]
\end{remark}

\begin{remark} Assumption \eqref{eq:conditions_on_V_1} can be replaced by $|V(0)|=M_0$. 
We require $V(0)$ to be positive
without loss of generality and only to simplify the notations in the proofs. 
\end{remark}

\begin{remark} Condition~\eqref{eq:conditions_on_V_4} may be relaxed if the forward and pullback attraction
hold for functions $v$ that satisfy certain localization conditions (belong to $L^1$, $L^2$, or have compact support).
\end{remark}

\subsection{Localized Gibbs distributions}
Let $\Gamma_{n_1,n_2}(x_1,x_2)$ be the set of all admissible paths on $[n_1,n_2]$ with fixed endpoints
$\gamma(n_1)=x_1,\ \gamma(n_2)=x_2$.

We say that a measure $\mu$ on $(\Z^d)^\Z$ (i.e. on paths $\alpha:\Z\to\Z^d$) is a Gibbs measure corresponding to
a realization of the potential $\phi:\Z\times\Z^d\to\R$ if it satisfies the DLR condition with 
Gibbs specification given by
\[
\mu_{n_1,n_2}(\gamma|\ x_1,x_2)=\frac{1}{Z_{n_1,n_2}(x_1,x_2)}e^{\Phi_{n_1,n_2}(\gamma)},\quad \gamma\in\Gamma_{n_1,n_2}(x_1,x_2), 
\]
with
\begin{equation*}
Z_{n_1,n_2}(x_1,x_2)=\sum_{\gamma\in\Gamma_{n_1,n_2}(x_1,x_2)}e^{\Phi_{n_1,n_2}(\gamma)}.
\end{equation*}
Namely, for any times $n_1<n_2$ and any points $x_1,x_2\in\Z$ 
the conditional distribution on paths in $\Gamma_{n_1,n_2}(x_1,x_2)$
defined by $\mu$ 
agrees with this specification: 
for any path $\gamma\in\Gamma_{n_1,n_2}(x_1,x_2)$, 
\begin{equation*}
\mu\left\{\alpha_{[n_1,n_2]}=\gamma \ |\ \alpha_{[n_1,n_2]}\in\Gamma_{n_1,n_2}(x_1,x_2)\right\}=
\mu_{n_1,n_2}(\gamma|\ x_1,x_2),
\end{equation*}
where $\alpha_{[n_1,n_2]}$ is the restriction of $\alpha$ on $[n_1,n_2]$.

Notice that $Z_{n_1,n_2}(x_1,x_2)$ can be expressed via the cocycle $T$:
\begin{equation}
Z_{n_1,n_2}(x_1,x_2)=(2d+1)^{n_2-n_1}T^{n_1,n_2}\delta_{x_1}(x_2),
\end{equation}
where
\begin{equation*}
\delta_{x_1}(x)=\begin{cases}1,&x=x_1,\\ 0,&x=0.\end{cases}
\end{equation*}


\begin{theorem}\label{th:gibbs-localiz} There is a set $\Omega'\subset \Omega$ of probability 1 with the following properties: 
\begin{enumerate}
\item  For every $\omega\in\Omega'$, there is a Gibbs measure $\mu=\mu^\omega$ corresponding to the realization
of the potential $\phi$. 
\item For every $\omega\in\Omega'$, the measure $\mu$ is a unique Gibbs measure with the following property: for every $\eps>0$
there is a number $r_\eps$ such that
\begin{equation}
\label{eq:weak_localization_condition}
 \liminf_{n\to\pm\infty} \mu\{\alpha:\ |\alpha_n|>r_\eps\}<\eps.
\end{equation}
\end{enumerate}
Moreover, for any $\lambda\in(0,\lambda_0)$, the measure $\mu$ is $2\lambda$-localized: there is a random variable $c=c(\omega)>0$
such that for every $\omega\in\Omega'$, every $n\in\Z$,
\begin{equation}
\label{eq:loc_for_gibbs}
 \mu\{\alpha: \alpha_n=x\}< c(\theta^n \omega) e^{-2\lambda |x|},\quad x\in\Z^d.
\end{equation}
\end{theorem}

\begin{remark}
 Notice 
that the exponent in the r.h.s.\ of~\eqref{eq:loc_for_gibbs} is $2\lambda$ rather than~$\lambda$ that
appears in the statement of Theorem~\ref{thm:perron-frobenius-cocycle}. This is due to two-sided estimates that we can use
in the proof of Theorem~\ref{th:gibbs-localiz}.
\end{remark}


A proof of Theorem~\ref{th:gibbs-localiz} is given in Section~\ref{sec:Gibbs}.

\section{Main localization lemma: the optimal path vs.\ entropy}\label{sec:main_loc_lemma}

This section is devoted to the estimate that plays the central role in the proofs given in forthcoming sections.
We shall need the following ``optimal'' path $\gamma^*$:
\begin{equation}
\label{eq:optimal_gamma}
\gamma^*_n(\omega)=
\begin{cases}0,& \omega_n=1\\
             e_1, &\omega_n=-1.
\end{cases}
\end{equation}
where $e_1=(1,0,\ldots,0)\in\Z^d$.
Clearly,
\[
 \Phi_{n_1,n_2}(\gamma^*)\ge \sum_{n_1+1}^{n_2}\xi_m,
\]
where the random variables $(\xi_m)_{m\in\Z}$ are defined by
\[
 \xi_m(\omega)=M_0\ONE_{\{\omega_m=1\}}-M_1\ONE_{\{\omega_m=-1\}}.
\]

\begin{lemma}\label{lm:main_localization_lemma}  Let $\lambda\in(0,\lambda_0)$. Then
there are random variables $\nu^-(\omega)=\nu^-(\ldots,\omega_{-1},\omega_0)$ and $\nu^+(\omega)=\nu^+(\omega_1,\omega_2,\ldots)$  such that with probability 1, for any $k\ge \nu^-(\omega)$,
\[
e^{\Phi_{-k,0}(\gamma^*)}\ge e^{\sum_{m=-k+1}^0\xi_m} > (2d+1)^k e^{M_1k} e^{\lambda k},
\]
and for any $k\ge \nu^+(\omega)$,
\[
e^{\Phi_{0,k}(\gamma^*)}\ge e^{\sum_{m=1}^{k}\xi_m} > (2d+1)^k  e^{M_1k} e^{\lambda k}.
\]
In addition, $\Pp\{\nu^-(\omega)=1\}>0$, and $\Pp\{\nu^+(\omega)=1\}>0$.  
\end{lemma}
\bpf Both  $\sum_{m=1}^{k}\xi_m$ and $\sum_{m=-k+1}^0\xi_m$ are
sums of i.i.d.\ random variables with mean equal to $(M_0-M_1)/2$. 
Let
\begin{equation}
\label{eq:epsilon}
\eps=\frac{M_0-M_1}{2}-M_1-\ln(2d+1)-\lambda.
\end{equation}
Notice that $\eps>0$ due to~\eqref{eq:conditions_on_V_3}. Now the strong law of large numbers implies the existence of random variables $\nu^\pm$
such that for almost every $\omega$ if $k\ge\nu^-(\omega)$ then 
\[
\sum_{m=-k+1}^0\xi_m >((M_0-M_1)/2-\eps)k.
\]
and if $k\ge\nu^+(\omega)$ then 
\[
\sum_{m=1}^{k}\xi_m >((M_0-M_1)/2-\eps)k.
\]
Therefore, for these values of $k$, respectively,
\begin{align*}
\frac{1}{(2d+1)^k} \cdot e^{\sum_{m=-k+1}^0\xi_m}>e^{((M_0-M_1)/2-\eps - \ln(2d+1))k}= e^{M_1k} e^{\lambda k},\\
\frac{1}{(2d+1)^k} \cdot e^{\sum_{m=1}^{k}\xi_m}>e^{((M_0-M_1)/2-\eps - \ln(2d+1))k}= e^{M_1k} e^{\lambda k},
\end{align*}
and the proof is complete.
\epf

\section{Proof of Theorem~\ref{thm:perron-frobenius-cocycle}}\label{sec:cocycle}

Our plan is  to find, for any large value of $r\in\N$, a sequence of numbers $(n_i)$  decreasing to $-\infty$ and such that, $T^{n_{i+1},n_i}$ is essentially a contraction in Hilbert projective metric on functions restricted to $B_r=[-r,r]^d$,
the ball of radius $r$ in the  sup-norm $|\cdot|$ on $\Z^d$.

The proof relies on several lemmas. The first part of the proof is devoted to finding an invariant and almost compact set 
essentially supporting the dynamics.

We begin with the main localization result. For any $c\ge 1$, we denote
\[
F(c)=\{\phi:\Z^d\to\R_+:\ \|\phi\|\le c\phi(0)\}. 
\]

\begin{lemma}\label{lm:radius_loc} Suppose $\lambda\in(0,\lambda_0)$. There is a number $K_1(\lambda)\ge 1$ such that if
$c\ge1$, $\phi\in F(c)$, $|y|\ge \nu^-(\omega)$, and $n\ge\nu^-(\omega)$, then 
\[
 T^{-n,0}\phi(y)\le K_1(\lambda) (e^{-\lambda|y|}+ce^{-\lambda n})T^{-n,0}\phi(0).
\]
\end{lemma}
 
The proof of this Lemma is given in Section~\ref{sec:auxi_lemmas}.

For $\lambda\in(0,\lambda_0)$, $r\in\N$, and $c\ge 1$, we define
\[
 n_0(\lambda,r,c)=\frac{\ln c}{\lambda}+r+1.
\]

\begin{lemma}\label{lm:initialization_in_the_past} Suppose $\lambda\in(0,\lambda_0)$, $r\in\N$, and $c\ge 1$.
If $n>n_0(\lambda,r,c)$,
$r\ge \nu^-(\omega)$, and $\phi\in F(c)$, then
\[
 \|(T^{-n,0}\phi)\ONE_{B_r^c}\|\le 2K_1(\lambda)e^{-\lambda (r+1)} T^{-n,0}\phi(0).
\]
\end{lemma}

This lemma is a direct implication of Lemma~\ref{lm:radius_loc} since its conditions automatically imply $n>\nu^-(\omega)$.

For the rest of the proof we fix $\lambda\in(0,\lambda_0)$. Lemma~\ref{lm:main_localization_lemma} implies that for any $r\in\N$, there is an event $A$ with $\Pp(A)>0$
such that for every $\omega\in A$ the following conditions hold true:

\begin{enumerate}
\item $\omega_{1}=\ldots=\omega_{r}=1$; \label{cnd:first_pluses-1.5}
\item $\nu_+(\theta^{r}\omega)=1$; \label{cnd:nu_plus_immediately-1.5}
\item $\omega_{-r+1}=\omega_{-r+2}=\ldots=\omega_{0}=1$; \label{cnd:terminal_pluses-1.5}
\item $\nu_-(\theta^{-r}\omega)=1$. \label{cnd:nu_minus_immediately-1.5}
\end{enumerate}

Therefore, with probability 1, we can choose a sequence $(n_i)_{i\in\N}$ (depending on $\omega$) decreasing to $-\infty$ and
such that the following Conditions 1--5 are satisfied for each $i\in\N$:
\begin{enumerate}
\item $\omega_{n_i+1}=\ldots=\omega_{n_i+r}=1$; \label{cnd:first_pluses-2}
\item $\nu_+(\theta^{n_i+r}\omega)=1$; \label{cnd:nu_plus_immediately-2}
\item $\omega_{n_i-r+1}=\omega_{n_i-r_0+2}=\ldots=\omega_{n_i}=1$; \label{cnd:terminal_pluses-2}
\item $\nu_-(\theta^{n_i-r}\omega)=1$. \label{cnd:nu_minus_immediately-2}
\item $n_{i-1}-n_{i}>2 n_0(\lambda,r,2K_1(\lambda))$.\label{cnd:spacing_in_sequence}
\end{enumerate}
The sequence can be chosen in a measurable way. Notice that if it satisfies conditions 1--5 for some $r$ it also satisfies the same conditions with $r$ replaced by
any nonnegative $r'<r$. 
 
For $r\in\N$, we define
\[
G(\lambda,r)=\{\phi\in F(2K_1(\lambda)):\ \|\phi\ONE_{B_r^c}\|\le \|\phi \ONE_{B_r}\|\}.
\]

\begin{lemma}\label{lm:G_preserved} There is a nonrandom $r_0\in\N$ such that  
\begin{enumerate}
 \item 
If $c\ge 1$,  $r>r_0$, and $-n<n_i-n_0(\lambda,r,c)$, then
\[
T^{-n,n_i} F(c)\subset G(\lambda,r).
\] 
\item If $r>r_0$, then for any $i\in\N$,
\[
T^{n_{i},n_{i-1}} G(\lambda,r)\subset G(\lambda,r). 
\]
\end{enumerate}
\end{lemma}

\bpf The first part of the lemma follows from Lemma~\ref{lm:initialization_in_the_past}.
The second part is a consequence of the first one and Condition~\ref{cnd:spacing_in_sequence}. \epf

\begin{lemma}\label{lm:boundedness_in_Hilbert} For $r\in\N$, there is a number $K_2(\lambda, r)\ge 1$ such that if
$\|\phi\ONE_{B_r^c}\|\le \|\phi\ONE_{B_r}\|$ then for any $y_1$ and $y_2$ with $|y_1|,|y_2|\le r$,
and any $i\in\N$,
\[
\frac{1}{K_2(\lambda, r)}\le\frac{T^{n_i,n_{i-1}}\phi(y_1)}{T^{n_i,n_{i-1}}\phi(y_2)}\le K_2(\lambda, r).
 \]
\end{lemma}

The proof is given in Section~\ref{sec:auxi_lemmas}.

We introduce now
\[
 H(\lambda,r)=\left\{\phi\in G(\lambda,r):\ \|\phi\|=1\ \text{and}\ \phi(y)\ge \frac{1}{K_2(\lambda,r)}\ \text{for}\ |y|\le r\right\}.
\]

\begin{lemma} \label{lm:image_F(c)subsetH}
\begin{enumerate}\item
If $r>r_0$, then for any $i\in\N$,
 \begin{align*} 
\wb T^{n_i,n_{i-1}} H(\lambda,r)\subset H(\lambda,r).
\end{align*}
\item\label{it:image_F(c)subsetH} If $c>0$,  $r>r_0$, $-n<n_i-n_0(\lambda,r,c)$, then
\[
 \wb T^{-n,n_{i-1}}F(c)\subset H(\lambda,r).
\]
\end{enumerate}
\end{lemma}
\bpf Lemmas~\ref{lm:G_preserved} and~\ref{lm:boundedness_in_Hilbert} imply 
\[
\wb T^{n_i,n_{i-1}} G(\lambda,r)\subset H(\lambda,r),
\]
and the first part of the lemma follows since $H(\lambda,r)\subset G(\lambda,r)$. The second part of the lemma follows from the first one and Lemma~\ref{lm:G_preserved}.\epf

The Hilbert projective metric $\rho_r$ between 
two functions $\phi,\psi:B_r\to\R_+$ is defined by
\begin{equation*}
\rho_r(\phi,\psi)=\ln\left(\max_{|x|\le r} \frac{\phi(x)}{\psi(x)}\cdot \max_{|x|\le r}\frac{\psi(x)}{\phi(x)}\right).
\end{equation*}
%
The Hilbert metric $\rho_r$ is not a true metric since it does not distinguish functions
proportional to each other. However it does define a metric on normalized positive functions defined on $B_r$.

For a function $\phi:\Z^d\to\Z_+$, we denote its restriction onto $B_r$ by $\pi_r\phi$. 

We notice that $\pi_r H(\lambda,r)$ is closed in $\rho_r$, and $\diam_r( H(\lambda,r))<\infty$, where for a set $A$ we denote
\[\diam_r (A)=\sup\{\rho_r(\pi_r\phi,\pi_r\psi):\ \phi,\psi\in A\}.\]

The following is the main contraction estimate:
\begin{lemma}
\label{lm:contraction_and_correction}  There are numbers $K_3(\lambda), K_4(\lambda)$ such that
for any $r$, any set $A\subset H(\lambda,r)$, and any $i\ge 2$,
\[
 \diam_r (T^{n_i,n_{i-1}} A)\le (1-K_3(\lambda)e^{-2\lambda_1 r})\diam_r(A)+K_4(\lambda)e^{-2\lambda_0r}.
\]
\end{lemma}

The proof is given in Section~\ref{sec:auxi_lemmas}. 

We can apply this estimate recursively along the sequence $(n_i)$. We need the following lemma to make connection to time 0 which is not included in~$(n_i)$:

\begin{lemma}\label{lm:terminal} There is a number $K_5(\lambda)$ such that  for any $r\in\N$ and
sufficiently large $i$, if $\phi,\psi\in F(2K_1(\lambda))$, 
then
\[
 \rho_r(\pi_rT^{n_i,0}\phi, \pi_r T^{n_i,0}\psi) \le \rho_r(\pi_r\phi,\pi_r\psi)+K_5(\lambda)e^{-\lambda_2 r},
\]
where
\[
\lambda_2=2(M_0-M_1-\ln(2d+1))>0.
\]
\end{lemma}

The proof is given in Section~\ref{sec:auxi_lemmas}.


\bpf[Proof of Theorem~\ref{thm:perron-frobenius-cocycle}] 
Lemma~\ref{lm:image_F(c)subsetH} implies that for sufficiently large $n$,
\begin{equation}
\label{eq:inclusionTH1}
 T^{-n,0} F(c)\subset T^{n_{i-1},0}H(\lambda,r),
\end{equation}
and if $i_1<i_2$,
\begin{equation}
\label{eq:inclusionTH2}
 T^{n_{i_2},0}H(\lambda,r)\subset T^{n_{i_1},0}H(\lambda,r).
\end{equation}

Take any sequence $\phi=(\phi_n)_{n\in\N}$ in $F(c)$ and consider the sequence $( T^{-n,0}\phi_n)$. 
Relations~\eqref{eq:inclusionTH1} and~\eqref{eq:inclusionTH2} imply that for each $r$, $\pi_r T^{-n,0}\phi_n$ are uniformly bounded
in $n$ and $\phi$. 
We consider the all pointwise limit points of $T^{-n,0}\phi_n$. Since pointwise limit points are not necessarily normalized in the uniform norm,
we normalize them and denote the resulting set by $S_\phi$. We denote by $S$ the union of $S_\phi$ over all possible
sequences $\phi$.   

It follows from the classical diagonal method that $S$ is not empty. Let us show that $S$, in fact, consists of a single point.

Suppose that on the contrary
$S$ contains at least two different points $\psi$ and $\psi'$. Then there exists $\bar \rho>0$ such that
$\rho_r(\pi_r\psi,\pi_r\psi')\geq \bar\rho$ for all large enough $r$. 
Since $\psi$ and $\psi'$ are
limiting points, there exist two sequences of functions $(\psi_{k})$ and $(\psi'_k)$ and two sequences of times
$(m_k)$ and $(m'_k)$ decaying to $-\infty$ such that for all large enough $k$,
\begin{equation}
\rho_r(\pi_rT^{-m_k,0}\psi_k,\pi_rT^{-m'_k,0}\psi'_{k})\geq \bar\rho/2. 
\label{eq:barrho}
\end{equation}
On the other hand, for any $i$, if $k$ is large enough, then both $T^{-m_k,n_i}\psi_k$ and $T^{-m'_{k},n_i}\psi'_k$
belong to $F(2K_1(\lambda))$. Lemmas~\ref{lm:contraction_and_correction} and~\ref{lm:terminal} imply that 
\begin{equation*}
\limsup_{i\to\infty} \diam_r(T^{-n_i,0}F(2K_1(\lambda)))\leq d(\lambda,r)=\frac{K_4(\lambda)e^{-2\lambda_0r}}{K_3(\lambda)e^{-2\lambda_1 r}}+K_5(\lambda)e^{-\lambda_2 r}.
\end{equation*} 
Taking $r$ large enough so that
$d(\lambda,r)<\bar\rho/2$, we obtain a contradiction with~\eqref{eq:barrho}. We conclude that $S_\infty$ cannot contain two distinct elements. Therefore, $S_\infty=\{\psi_\infty\}$ for some $\phi_\infty$, and it is easy to see
that $\psi_\infty$ does not depend on $c$.

We now set $u^\omega=\psi_\infty$, where $\psi_\infty$ is the unique element of $S$. The uniqueness
above ensures that $u$ is a positive cocycle eigenfunction satisfying~\eqref{eq:uniqueness_condition}. It is also obviously unique. The desired localization property follows from Lemma~\ref{lm:radius_loc}.
The pullback attraction follows since for any $\eps$, we can find $r_0$ and $n_0$ such that
\[
 \|u\ONE_{B_{r_0}^c}\|< \eps,
\]
\[
 \|(\wb T^{-n,0}v)\ONE_{B_{r_0}^c}\|< \eps,\quad n>n_0,
\]
and
\[
 \rho_{r_0}(\wb T^{-n,0} v,u)<\eps, \quad n>n_0.
\]
The forward attraction is proven similarly.\epf

\section{Proof of Theorem~\ref{th:gibbs-localiz}}\label{sec:Gibbs}
For a fixed $\omega$, we say that a measure $\mu_{n_1,n_2}$ on $\Gamma_{n_1,n_2}$ is a finite volume Gibbs distribution on $[n_1,n_2]$ for the realization of the potential $\phi$ if
for any points $x_1,x_2$ and any path $\gamma\in\Gamma_{n_1,n_2}(x_1,x_2)$,
\[
 \mu_{n_1,n_2}\bigl(\{\gamma\}\ |\ \Gamma_{n_1,n_2}(x_1,x_2)\bigr)=\frac{e^{\Phi_{n_1,n_2}(\gamma)}}{Z_{n_1,n_2}(x_1,x_2)}.
\]

Let us introduce $\nu=\max\{\nu^+,\nu^-\}$, where $\nu^\pm$ are introduced in Section~\ref{sec:main_loc_lemma}.

The proof of the following localization lemma for finite volume Gibbs distributions is given in Section~\ref{sec:auxi_lemmas}:
\begin{lemma}\label{lm:concentration_for_Gibbs_marginals} For any $\lambda\in(0,\lambda_0)$, there is a constant $C>0$ and a function $N_1:\N\to\N$ with the following property.
Suppose $\mu_{n_1,n_2}$ is a finite volume Gibbs measure on an interval $[n_1,n_2]$ for a realization of the potential $\phi$. If 
$r\in\N$ and $n\in[n_1,n_2]$ satisfy \[\mu_{n_1,n_2}\{\alpha:\ \alpha_{n_1},\alpha_{n_2}\in B_r\}=1,\] 
\[
\nu(\theta^n\omega)<r,
\]
\[
n_2-n>N_1(r),\]\[n-n_1>N_1(r),
\]
then
\[
\mu_{n_1,n_2}\{\alpha:\ |\alpha_n|>r\}< C e^{-2 \lambda r}.
\]
\end{lemma}

\medskip

Let us prove the existence first. 
We fix an $\omega\in\Omega$ and for each $m\in\N$ consider a unique measure $\mu^m$ on $(\Z^d)^\Z$ such that 
\begin{enumerate}
\item $\mu^m\{\alpha:\ \alpha_k=0\}=1$ for all $k\le -m$ and all $k\ge m$ 
\item The projection of $\mu^m$ on $[-m,m]$ is a finite volume Gibbs measure.
\end{enumerate}

The following statement is a direct consequence of Lemma~\ref{lm:concentration_for_Gibbs_marginals}.

\begin{lemma}\label{lm:basic_localization_lemma2} For any $\lambda\in(0,\lambda_0)$, there is a random variable $c(\omega)$ 
such that for almost every $\omega\in\Omega$ and every $r$,
\begin{equation*}
\mu^m\{\alpha: |\alpha_n|>r\}\le c(\theta^n \omega)e^{-2\lambda r}.
\end{equation*} 
\end{lemma}

Applying this result, we conclude that with probability 1, the sequence of measures $(\mu^m)$ is tight in $(\Z^d)^{\Z}$, and, due to the Prokhorov criterion, it contains a weakly convergent subsequence. Denoting the limit of this subsequence by $\mu$, one can easily verify that
$\mu$ is a Gibbs measure satisfying~\eqref{eq:loc_for_gibbs}.

\medskip

To prove the uniqueness, we must show that any two Gibbs measures $\mu$ and $\mu'$ satisfying~\eqref{eq:weak_localization_condition} coincide. The plan is as follows. We shall consider a sequence of restrictions of $\mu$ and $\mu'$ on intervals $[n_j^-,n_j^+]$, with $n_j^+\to\infty$
and
$n_j^-\to -\infty$ as $j\to\infty$. We shall iteratively estimate the proximity of these restrictions to each other in total variation, by showing
 that the restrictions on  $[n_{j-1}^-,n_{j-1}^+]$  are  (up to a small correction) closer to each other than the respective restrictions on $[n_{j}^-,n_{j}^+]$,
by a multiplicative factor that is less than~1. The multiplier and the correction can be controlled by the
choice of the sequences $(n_j^+)$ and $(n_j^-)$.

\begin{lemma}\label{lm:coupling_minorant} There is a constant $c>0$ and a function $N_2:\N\to\N$ with the following property. Suppose $r\in\N$, , and
\[
\omega_{-N_2(r)+1}=\omega_{-N_2(r)+2}=\ldots=\omega_{N_2(r)}=1.
\] If $n>N_2(r)$ and $\mu$ a finite volume Gibbs measure on $[-n,n]$ such that
\[
\mu\{\gamma:\ \gamma_{-N_2(r)},\gamma_{N_2(r)}\in B_r\}=1, 
\]
 then
\[
 \mu\{\gamma:\ \gamma_0=0\}>c \mu\{\gamma:\ \gamma_0=x\},\quad x\ne 0. 
\]
\end{lemma}
The proof is analogous to that of Lemma~\ref{lm:radius_loc}.
\begin{lemma}\label{lm:coupling_times} Let $\lambda\in(0,\lambda_0)$. Then for almost every $\omega\in\Omega$ and for every $r>0$,
there a doubly infinite sequence $(n_i)_{i\in\Z}$ such that 
\[\lim_{i\to+\infty} n_i=\infty,\quad \lim_{i\to-\infty} n_i=-\infty,\]  
and for any $i$, 
\[
\nu^-(\theta^{-N_2(r)+n_i}\omega)=1,
\]
\[
\omega_{-N_2(r)+1+n_i}=\omega_{-N_2(r)+2+n_i}=\ldots=\omega_{N_2(r)+n_i}=1,
\]
\[
 \nu^+(\theta^{N_2(r)+n_i}\omega)=1,
\]
and for all $i$, 
\[n_{i+1}-n_i>2N_2(r).\]
\end{lemma}
\bpf Given $r>0$.
\begin{align*}
 \Pp\bigl\{\nu^-(\theta^{-N_2(r)}\omega)=1;\ \omega_{-N_2(r)+1}=\ldots=\omega_{N_2(r)}=1;\nu^+(\theta^{N_2(r)}\omega)=1\bigr\}
\\=\Pp\{\nu^-(\theta^{-N_2(r)}\omega)=1\}\Pp\{\omega_{-N_2(r)+1}=\omega_{-N_2(r)+2}=\ldots=\omega_{N_2(r)}=1\}\\
\times \Pp\{\nu^+(\theta^{N_2(r)}\omega)=1\}>0,
\end{align*}
so that the lemma follows from the Bernoulli property.
\epf

Now we return to the proof of the uniqueness in Theorem~\ref{th:gibbs-localiz}.
For a measure $\mu$ on $(\Z^d)^\Z$ and any set $S\subset\Z$ we denote by $\mu_S$  the measure induced by $\mu$ on paths restricted to $S$.

Let now $\mu$ and $\mu'$ be two Gibbs measures satisfying~\eqref{eq:weak_localization_condition} for a given $\omega$.

We have to show that for any $l>0$, the distributions induced by $\mu$ and $\mu'$ on trajectories defined on $[-l,l]$ coincide: $\mu_{[-l,l]}=\mu'_{[-l,l]}$. Since $\mu$ and $\mu'$ are Gibbs measures with nearest neighbor interaction,  it is sufficient to
check that the two-dimensional boundary distributions coincide: $\mu_{\{-l,l\}}=\mu'_{\{-l,l\}}$.

We fix an arbitrary $\eps>0$ and use~\eqref{eq:weak_localization_condition} to find $r_\eps>0$ and sequences $(m_k)_{k\in\Z}$, $(m'_k)_{k\in\Z}$ 
such that
\[
 \mu\bigl\{\alpha:\ |\alpha_{m_k}|>r_\eps\bigr\}<\eps,\quad k\in\Z,
\]
and
\[
 \mu'\bigl\{\alpha:\ |\alpha_{m_k}|>r_\eps\bigr\}<\eps,\quad k\in\Z.
\]
For $r\ge r_\eps$ and $k\in\N$ consider measures $\mu^{k,r}$ and $\mu'^{k,r}$ obtained from
$\mu$ and $\mu'$ by conditioning, respectively, on $\{\alpha_{m^\pm_k}\le r\}$ and $\{\alpha_{m'^\pm_k}\le r\}$. 
Due to the arbitrary choice of $\eps$, it is sufficient to show that the total variation distance
between $\mu^{k,r}_{\{-l,l\}}$ and $\mu'^{k,r}_{\{-l,l\}}$ can be made arbitrarily small by choosing sufficiently large $r$ and $k$.

Let us fix $r$ and find the sequence $(n_i)_{i\in\Z}$ provided by Lemma~\ref{lm:coupling_times}. For any given $i_0\in\N$, one can find
$k$ such that \[n_{i_0}+N(r)<\min\{m_k,m'_k\}\] and \[n_{-i_0}-N(r)>\max\{m_{-k},m'_{-k}\}.\]   
Lemma~\ref{lm:concentration_for_Gibbs_marginals} implies that if $|i|<i_0$, then
\begin{align}
\mu^{k,r}\{\alpha:\ |\alpha_{n_i}|>r\}< C e^{-2 \lambda r},\label{eq:loc_gibbs1}\\
\mu^{k,r}\{\alpha:\ |\alpha_{n_i+N_2(r)}|>r\}< C e^{-2 \lambda r},\label{eq:loc_gibbs2}\\
\mu^{k,r}\{\alpha:\ |\alpha_{n_i-N_2(r)}|>r\}< C e^{-2 \lambda r},\label{eq:loc_gibbs3}
\end{align}
and same estimates hold for $\mu'^{k,r}$.
For $0<i<i_0$, let us denote
\[
 \mu_i=\mu_{\{n_{-i},n_{i}\}}^{k,r},\quad \mu'_i=\mu'^{k,r}_{\{n_{-i},n_{i}\}}.
\]
We are going to estimate the total variation distance $d_{TV}(\mu_{i-1},\mu'_{i-1})$ via $d_{TV}(\mu_{i},\mu'_{i})$.

For any $i$ we introduce $\sigma_i$ to be the maximal measure minorizing both $\mu_i$ and $\mu'_i$ and concentrated on $B_r\times B_r$: 
\begin{align*}
\mu_i&=\sigma_i+\rho_i+\delta_i\\
\mu'_i&=\sigma_i+\rho'_i+\delta'_i,
\end{align*}
Here $\rho_i$ and $\rho'_i$ are mutually singular measures on $B_r\times B_r$, and measures $\delta_i,\delta'_i$ are supported on
$\Z^d\times \Z^d\setminus B_r\times B_r$.
We have
\[
d_{TV}(\mu_{i},\mu'_{i})\ge\rho_i(B_r),
\]
and, due to~\eqref{eq:loc_gibbs2},\eqref{eq:loc_gibbs3},
\begin{align*}
 \delta_i(\Z^d\times\Z^d)\le 2Ce^{-2\lambda r},\\
\delta'_i(\Z^d\times\Z^d)\le 2Ce^{-2\lambda r}.
\end{align*}

Combining a basic coupling estimate based on Lemma~\ref{lm:coupling_minorant}, the fact that there are $(2r+1)^{2}$ points in $B_r$,
and the estimate~\eqref{eq:loc_gibbs1},
we obtain:
\[
d_{TV}(\mu_{i-1},\mu'_{i-1})\le d_{TV}(\mu_{i},\mu'_{i})\left(1-\left(\frac{c}{(2r+1)^d-1+c}\right)^2\right)+2Ce^{-2\lambda r}.
\]

Since the total variation distance is always bounded by 1,  applying this estimate iteratively, we obtain
\[
d_{TV}(\mu_i,\mu'_i)< 2\frac{C}{c^2}e^{-2\lambda r}((2r+1)^d-1+c)^2,\quad i\in\N.
\]
Since the choice of $r$ is arbitrary and the r.h.s.\ converges to 0 as $r\to\infty$, the proof is completed.\epf

\section{Proofs of auxiliary lemmas}\label{sec:auxi_lemmas}
\bpf[Proof of Lemma \ref{lm:radius_loc}]
Let us write 
\begin{align}\notag
(2d+1)^{n}T^{-n,0}\phi(y)&\le \sum_{k=|y|}^{n}\sum_x \phi(x) Z_{-n,-k}(x,0) \tilde Z_{-k,0}(0,y)+\bar Z_{-n,0}(y)\|\phi\|\\
&\le \sum_{k=|y|}^{n} I_k +\bar I.\label{eq:decomposition_in_I_k}
\end{align}
Here $\tilde Z_{-k,0}(0,y)$ is the partition function over all paths on $[-k,0]$ connecting $0$ to $y$ and avoiding $0$ after time $-k$;
$\bar Z_{-n,0}(y)$ is the partition function over paths on $[-n,0]$ ending up at $y$ and never visiting $0$.
 
Considering all paths that coincide with $\gamma^*$ on $[-k+1,0]$, we can write 
\[
(2d+1)^{n}T^{-n,0}\phi(0)\ge  \sum_{x}\phi(x)Z_{-n,-k}(x,0)e^{\sum_{m=-k+1}^{0}\xi_m}e^{-(M_0+M_1)},
\]
where the factor of $e^{-(M_0+M_1)}$ appears since it is possible that  $\gamma^*(0)\ne 0$. 
We also notice that
\[
\tilde Z_{-k,0}(0,y)\le (2d+1)^{k}e^{M_1k} 
\]
since there are at most $(2d+1)^{k}$ paths contributing to $\tilde Z_{-k,0}(0,y)$ and each path contributes at most $e^{M_1k}$. 
Combining these two estimates we get
\begin{align}\notag
 I_k&\le e^{M_0+M_1}\frac{(2d+1)^{k}e^{M_1k}}{e^{\sum_{m=-k+1}^{0}\xi_m}}(2d+1)^nT^{-n,0}\phi(0)\\
&\le e^{M_0+M_1} e^{-\lambda k}(2d+1)^nT^{-n,0}\phi(0),
\label{eq:I_k-improved}
\end{align}
where the second inequality follows from Lemma~\ref{lm:main_localization_lemma} and condition $k\ge|y|\ge \nu^-(\omega)$.
Similarly, for $\bar I$ we get:
\begin{align*}
\bar I=\bar Z_{-n,0}(y)\|\phi\|\le (2d+1)^n e^{M_1n}\|\phi\|\le (2d+1)^n e^{M_1n}c\phi(0). 
\end{align*}
On the other hand,
\begin{align*}
(2d+1)^n T^{-n,0}\phi(0)&\ge \phi(0) e^{\sum_{m=-n+1}^{0}\xi_m}e^{-(M_0+M_1)}
\\ &\ge \phi(0) (2d+1)^{n}e^{M_1n}e^{\lambda n}e^{-(M_0+M_1)},
\end{align*}
where the second inequality follows from Lemma~\ref{lm:main_localization_lemma} and condition $n\ge \nu^-(\omega)$.
Combining these estimates, we get
\begin{align}\notag
 \bar I
&\le (2d+1)^n e^{M_1n}c \frac{(2d+1)^n T^{-n,0}\phi(0)}{(2d+1)^{n}e^{M_1n}e^{\lambda n}e^{-(M_0+M_1)}}\\
&\le c e^{M_0+M_1} e^{-\lambda n}(2d+1)^n T^{-n,0}\phi(0).\label{eq:barI}
\end{align}
The lemma now follows from~\eqref{eq:decomposition_in_I_k},\eqref{eq:I_k-improved},and \eqref{eq:barI}.
\epf

\bpf[Proof of Lemma~\ref{lm:boundedness_in_Hilbert}] 
It is sufficient to prove the upper bound, so we write
\begin{align}\notag
(2d+1)^{n_{i-1}-n_i} &T^{n_i,n_{i-1}}\phi(y_1)\\ \le &\left(\sum_{n_i+r+1\le k\le l\le n_{i-1}-r} Z_{n_i,n_{i-1}}^{k,l}(y_1) + \bar Z_{n_i,n_{i-1}}(y_1)\right)\|\phi\|,
\label{eq:decomp_Z_again-for-boundedness-in-hilbert-norm}
\end{align}
where $Z_{n_i,n_{i-1}}^{k,l}(y_1)$ is the partition function over all paths $\gamma$ defined on interval~$[n_i,n_{i-1}]$ that terminate at~$y_1$ and satisfy
\[\min\{m> n_i+r: \gamma(m)=0\}=k,\]
\[\max\{m\le n_{i-1}-r: \gamma(m)=0\}=l,\]
and $\bar Z_{n_i,n_{i-1}}(y_1)$ is the partition function over all paths that terminate at~$y_1$ and do not visit 0 between $r$
and $n-r$.
We have
\begin{align*}
\notag
Z_{n_i,n_{i-1}}^{k,l}(y_1)\le& (2d+1)^{r}e^{M_0r}(2d+1)^{k-(n_i+r)} e^{M_1(k-(n_i+r)-1)}e^{M_0}\\&\times Z_{k,l}(0,0)(2d+1)^{n_{i-1}-r-l}e^{M_1(n_{i-1}-r-l)} (2d+1)^{r}e^{M_0r},
\end{align*}
and, considering a point $x^*$ such that $|x^*|\le r$ and $\phi(x^*)=\|\phi\|$,
\begin{align}
\notag
(2d&+1)^{n_{i-1}-n_i} T^{n_i,n_{i-1}}\phi(y_2)\\ &\ge \notag e^{-M_0r}e^{\sum_{m=n_i+r+1}^{k-1}\xi_m}e^{-M_0}Z_{k,l}(0,0)e^{\sum_{m=l+1}^{n_{i-1}-r}\xi_m}e^{-M_0r}\phi(x^*),
\end{align}
so that
\begin{multline}
\frac{Z_{n_i,n_{i-1}}^{k,l}(y)\|\phi\|}{(2d+1)^{n_i-n_{i-1}} T^{n_i,n_{i-1}}\phi(y_2)}\\ \le (2d+1)^{2r+2}e^{M_0(4r+2)}e^{-\lambda(k-(n_i+r)-1)} e^{-\lambda(n-r-l)}. \label{eq:Z_four_indices_lower_bound}
\end{multline}
Since
\[
 \sum_{n_i+r+1\le k\le l\le n_{i-1}-r}e^{-\lambda(k-(n_i+r)-1)} e^{-\lambda(n-r-l)}<\infty,
\]
the lemma follows from inequality~\eqref{eq:Z_four_indices_lower_bound} combined with an analogous estimate for $\bar Z_{n_i,n_{i-1}}(y_1)$ 
and~\eqref{eq:decomp_Z_again-for-boundedness-in-hilbert-norm}.
\epf

Now we begin preparations for the proof of the main contraction estimate, Lemma~\ref{lm:contraction_and_correction}.
For any $i$ and $r$, we define truncated operators
\begin{align*}
T_r^{n_i,n_{i-1}} u(y)&=\frac{1}{(2d+1)^{n_{i-1}-n_i}}\sum_{|x|\le r}Z_{n_i,n_{i-1}}(x,y)u(x)=T^{n_i,n_{i-1}}(u\ONE_{B_r})(y),\\
\hat T^{n_i,n_{i-1}} u(y)&=\frac{1}{(2d+1)^{n_{i-1}-n_i}}\sum_{|x|> r}Z_{n_i,n_{i-1}}(x,y)u(x)=T^{n_i,n_{i-1}}(u\ONE_{B_r^c})(y).
\end{align*}
The restriction of $T^{n_i,n_{i-1}}$ on $B_r$ can be viewed as a linear finite-dimensional operator in $\R^{B_r}$
given by a matrix   $(2d+1)^{n_{i}-n_{i-1}}Z_{n_i,n_{i-1}}(x,y)$, $x,y\in B_r$, with positive entries.
Therefore, we can apply a classical estimate on contraction in Hilbert metric, see e.g.,
Theorem 3.12 in \cite{SenetaMR719544}: for any functions $\phi,\psi:\Z^d\to\R_+$,
\begin{equation}
\rho_r(T_r^{n_i,n_{i-1}}\phi,T_r^{n_i,n_{i-1}}\psi) \le \frac{1-\sqrt{L^{n_i,n_{i-1}}_r}}{1+\sqrt{L^{n_i,n_{i-1}}_r}}\rho_{r}(\phi,\psi),
\label{eq:contraction_via_L}
\end{equation}
where
\begin{equation*}
L^{n_i,n_{i-1}}_r=\min_{|x_1|,|x_2|,|y_1|,|y_2|\le r}\left(\frac{Z_{n_i,n_{i-1}}(x_1,y_1)}{Z_{n_i,n_{i-1}}(x_2,y_1)}\cdot 
\frac{Z_{n_i,n_{i-1}}(x_2,y_2)}{Z_{n_i,n_{i-1}}(x_1,y_2)}\right).
\end{equation*}

\begin{lemma}\label{lm:contraction_coef} There is a positive constant $K_6(\lambda)$ such that
for any $i$ and $r$,
\[
L^{n_i,n_{i-1}}_r\ge K_6(\lambda) e^{-4\lambda_1 r}.
\] 
\end{lemma}

Since the dynamics of the actual system is not restricted to $B_r$ we have to estimate the influence of $B_r^c$:
\begin{lemma}\label{lm:influx} There is a positive number $K_7(\lambda)$ such that if 
$|y|\le r$, then
\[
 \sum_{x:|x|>r} Z_{n_i,n_{i-1}}(x,y)\le K_7(\lambda)e^{-2\lambda_0 r} Z_{n_i,n_{i-1}}(0,y).
\]
\end{lemma}

We postpone the proof of  Lemmas~\ref{lm:contraction_coef} and~\ref{lm:influx}  till the end of this section.

\bpf[Proof of Lemma~\ref{lm:contraction_and_correction}] Lemma~\ref{lm:contraction_coef} and estimate \eqref{eq:contraction_via_L} imply
that for any positive functions $\phi$ and $\psi$,
\begin{equation}
\label{eq:contraction_constants}
\rho_r(T_r^{n_i,n_{i-1}}\phi,T_r^{n_i,n_{i-1}}\psi)<\left(1-\sqrt{K_6(\lambda)}e^{-2\lambda_1}\right) \rho_r(\phi,\psi).
\end{equation}
To estimate the full untruncated operators, we write:
\begin{align}\label{eq:hilbert-metric}
\rho_r(T^{n_i,n_{i-1}}\phi,T^{n_i,n_{i-1}}\psi)=&
\ln\left(\max_{|y|\le r} \frac{T^{n_i,n_{i-1}}\phi(y)}{T^{n_i,n_{i-1}}\psi(y)}\cdot 
\max_{|y|\le r}\frac{T^{n_i,n_{i-1}}\psi(y)}{T^{n_i,n_{i-1}}\phi(y)}\right)
\\\notag \le& 
\ln\left(\max_{|y|\le r} \frac{T_r^{n_i,n_{i-1}}\phi(y)+\hat T_r^{n_i,n_{i-1}}\phi(y)}{T_r^{n_i,n_{i-1}}\psi(y)}\right)
\\ \notag &+\ln\left(
\max_{|y|\le r}\frac{T_r^{n_i,n_{i-1}}\psi(y)+\hat T_r^{n_i,n_{i-1}}\psi(y)}{T_r^{n_i,n_{i-1}}\phi(y)}\right)
\end{align}
Lemma \ref{lm:influx} implies
\begin{align*}
\hat T_r^{n_i,n_{i-1}}\phi(y)&\le \frac{1}{(2d+1)^{n_{i-1}-n_i}}K_7(\lambda)e^{-2\lambda_0 r} Z_{n_i,n_{i-1}}(0,y)\|\phi\ONE_{B_r}\|\\
&\le  K_7(\lambda)e^{-2\lambda_0 r} \frac{1}{(2d+1)^{n_{i-1}-n_i}} Z_{n_i,n_{i-1}}(0,y)\cdot 2K_1(\lambda)\phi(0)\\
&\le 2K_1(\lambda) K_7(\lambda)e^{-2\lambda_0 r} T_r^{n_i,n_{i-1}}\phi(y).
\end{align*}
Analogously,
\begin{align*}
\hat T_r^{n_i,n_{i-1}}\psi(y)&\le 2K_1(\lambda) K_7(\lambda)e^{-2\lambda_0 r} T_r^{n_i,n_{i-1}}\psi(y).
\end{align*}
Plugging the last two inequalities into~\eqref{eq:hilbert-metric}, we get 
\begin{align*}
\rho_r(T_r^{n_i,n_{i-1}}\phi,T_r^{n_i,n_{i-1}}\psi) &\le 
\ln\left(\max_{|y|\le r} \frac{T_r^{n_i,n_{i-1}}\phi(y)(1+2K_1(\lambda)K_7(\lambda)e^{-2\lambda_0 r})}{T_r^{n_i,n_{i-1}}\psi(y)}\right) 
\\&+\ln\left(\max_{|y|\le r}\frac{T_r^{n_i,n_{i-1}}\psi(y)(1+2K_1(\lambda)K_7(\lambda)e^{-2\lambda_0 r})}{T_r^{n_i,n_{i-1}}\phi(y)}\right)\\ 
&\le \rho_r(T_r^{n_i,n_{i-1}}\phi,T_r^{n_i,n_{i-1}}\psi)+ 4K_1(\lambda)K_7(\lambda)e^{-2\lambda_0 r}
\end{align*}
This estimate along with~\eqref{eq:contraction_constants} implies the lemma with $K_3(\lambda)=\sqrt{K_6(\lambda)}$ and $K_4(\lambda)=4K_1(\lambda)K_7(\lambda)$.\epf

\bpf[Proof of Lemma~\ref{lm:contraction_coef}]
Let us estimate the ratios in the r.h.s.\ of the definition of $L^{n_i,n_{i-1}}_r$.
For $n_i+r< k\le l\le n-r$, we define $Z_{n_i,n_{i-1}}^{k,l}(x,y)$ as the partition function over all paths $\gamma$ on $[n_i,n_{i-1}]$ connecting $x$ to $y$
satisfying 
\[\min\{m> n_i+r:\ \gamma(m)=0\}=k,\]
\[\max\{m\le n_{i-1}-r:\ \gamma(m)=0\}=l,\]
we also define $\bar Z_{n_i,n_{i-1}}(x,y)$ to be the partition function over all paths on $[n_i,n_{i-1}]$ connecting $x$ to $y$ that do not visit~0 between $n_{i}+r+1$ and $n_{i-1}-r$, so that
\begin{equation}\label{eq:decompositionwrttoching_zero}
Z_{n_i,n_{i-1}}(x,y)=\sum_{n_i+r < k \le l \le n_{i-1}-r} Z_{n_i,n_{i-1}}^{k,l}(x,y) + \bar Z_{n_i,n_{i-1}}(x,y).
\end{equation}
For $n_i+r< k\le l\le n_{i-1}-r$ and any points $x,y$ with $|x|,|y|\le r$,
\begin{align}\notag
Z_{n_i,n_{i-1}}^{k,l}(x,y)\le &(2d+1)^{|x|-1}e^{M_1(|x|-1)}(2d+1)^{r-|x|+1}e^{M_0(r-|x|+1)}\\ \notag
&\times (2d+1)^{k-1-(n_i+r)}e^{M_1(k-1-(n_i+r))}e^{M_0}Z_{k,l}(0,0)\\\notag&\times (2d+1)^{n-r-l}e^{M_1(n-r-l)}\\
&\times (2d+1)^{r-|y|}e^{M_0(r-|y|)}(2d+1)^{|y|}e^{M_1|y|}.\label{eq:sum_for_touching_0}
\end{align}

On the other hand, considering all paths that start at $x$, go straight to~$0$ (which takes $|x|$ steps),
stay at $0$ until time $n_i+r$ accumulating $M_0$ at each time step, follow the optimal path $\gamma^*$ from
$n_i+r+1$ to $k-1$ (if $k=n_i+r+1$ this part is empty), at $k$ visit $0$, return to~$0$ at time $l$, follow $\gamma^*$ up to $n_{i-1}-r$, stay at $0$ up to
$n_{i-1}-|y|$, go straight to $y$ where they terminate at time $n_{i-1}$, 
\begin{align}
\notag Z_{n_i,n_{i-1}}(x,y)\ge &e^{-M_1(|x|-1)}e^{M_0(r-|x|+1)}e^{\sum_{m=n_i+r+1}^{k-1} \xi_m} e^{-M_0}\\ &\times Z_{k,l}(0,0)
e^{\sum_{m=l+1}^{n_{i-1}-r}\xi_m} e^{M_0(r-|y|)}e^{-M_1|y|}.\label{eq:partition_along_0}
\end{align}
The definitions of $\nu^+$ and $\nu^-$ from Lemma~\ref{lm:main_localization_lemma} imply that,
due to our assumptions on $n_i$ and $n_{i-1}$,
\begin{equation}
\label{eq:lln-consequence-1}
 e^{\sum_{m=n_i+r+1}^{k-1} \xi_m}\ge (2d+1)^{k-1-(n_i+r)}e^{M_1(k-1-(n_i+r))} e^{\lambda (k-1-(n_i+r))},
\end{equation}
and
\begin{equation}
\label{eq:lln-consequence-2}
 e^{\sum_{m=l+1}^{n_{i-1}-r} \xi_m}\ge (2d+1)^{n_{i-1}-r-l}e^{M_1(n_{i-1}-r-l)}e^{\lambda(n_{i-1}-r-l)},
\end{equation}
so that
\begin{align}\notag
\frac{Z_{n_i,n_{i-1}}^{k,l}(x_2,y_1)}{Z_{0,n}(x_1,y_1)}
\le& (2d+1)^{2r+1}e^{M_1(|x_1|+|x_2|)}e^{M_0(|x_1|-|x_2|)}e^{M_0+M_1}
\\&\times e^{2 M_1|y_1|}e^{-\lambda((k-r)+(n-r-l)-2)}\notag
\\
\le& (2d+1)^{2r+1}e^{4M_1r}e^{M_0(|x_1|-|x_2|+2)} \notag\\  &\times e^{-\lambda((k-1-(n_i-r))+(n_{i-1}-r-l))}. \label{eq:generic_term_fraction}
\end{align}
For $\bar Z_{n_i,n_{i-1}}(x,y)$ we have
\begin{align*}
\bar Z_{n_i,n_{i-1}}(x,y)\le &(2d+1)^{|x|-1}e^{M_1(|x|-1)}(2d+1)^{r-|x|+1}e^{M_0(r-|x|+1)}\\
&\times (2d+1)^{n_{i-1}-n_{i}-2r}e^{M_1(n_{i-1}-n_{i}-2r)} \\
&\times (2d+1)^{r-|y|}e^{M_0(r-|y|)}(2d+1)^{|y|}e^{M_1|y|},
\end{align*}
and
\begin{align*}
Z_{n_i,n_{i-1}}(x,y)\ge &e^{-M_1(|x|-1)}e^{M_0(r-|x|+1)}e^{\sum_{m=n_i+r+1}^{n_{i-1}-r} \xi_m} e^{M_0(r-|y|)}e^{-M_1|y|},
\end{align*}
so that
\begin{align}
\frac{\bar Z_{n_i,n_{i-1}}(x_2,y_1)}{Z_{n_i,n_{i-1}}(x_1,y_1)}\le &(2d+1)^{2r+1}e^{4M_1r}e^{M_0(|x_1|-|x_2|)}.\label{eq:terminal_term_fraction}
\end{align}
Plugging \eqref{eq:generic_term_fraction} and \eqref{eq:terminal_term_fraction} into~\eqref{eq:decompositionwrttoching_zero}, we see that
\begin{align*}
\frac{Z_{n_i,n_{i-1}}(x_2,y_1)}{Z_{n_i,n_{i-1}}(x_1,y_1)}&\le (2d+1)^{2r+1}e^{4M_1r}e^{M_0(|x_1|-|x_2|)}K_8(\lambda),
\end{align*}
for some $K_8(\lambda)>0$. This inequality and its counterpart with $x_1,x_2$, and, respectively, $y_1$ replaced by $x_2,x_1$, and, respectively, $y_2$, immediately
implies:
\[
L^{n_i,n_{i-1}}_r\ge (2d+1)^{-4r-2}e^{-8M_1r} K^{-2}_8(\lambda),
\]
and the lemma holds true with $K_6(\lambda)=(2d+1)^{-2}K^{-2}_8(\lambda)$.
\epf

\bpf[Proof of Lemma~\ref{lm:influx}] The following decomposition is analogous to
~\eqref{eq:decompositionwrttoching_zero}:
\begin{equation}
\label{eq:decompositionwrttoching zero2}
\sum_{x:|x|>r} Z_{n_i,n_{i-1}}(x,y)\le \sum_{n_i+r<k\le l\le n_{i-1}-r} Z_{n_i,n_{i-1}}^{k,l}(y)+ \bar Z_{n_i,n_{i-1}}(y),  
\end{equation}
where for $n_i+r<k\le l\le n_{i-1}-r$,  $Z_{n_i,n_{i-1}}^{k,l}(y)$ is the partition function over all paths $\gamma$ on $[n_i,n_{i-1}]$ terminating at $y$ and satisfying 
\[\min\{m:\ \gamma(m)=0\}=k,\]
\[\max\{m\le n_{i-1}-r:\ \gamma(m)=0\}=l,\]
and $\bar Z_{n_i,n_{i-1}}(y)$ is the partition function over all paths terminating at $y$ that do not visit 0 between $n_i$ and $n_{i-1}-r$.
Analogously  to \eqref{eq:sum_for_touching_0},
\begin{align}
Z_{n_i,n_{i-1}}^{k,l}(y)\le &(2d+1)^{r}e^{M_1r} (2d+1)^{k-1-(n_i+r)}e^{M_1(k-1-(n_i+r))}e^{M_0}\notag
\\&\times Z_{k,l}(0,0) (2d+1)^{n_{i-1}-r-l}e^{M_1(n_{i-1}-r-l)}\notag\\
&\times (2d+1)^{r-|y|}e^{M_0(r-|y|)}(2d+1)^{|y|}e^{M_1|y|}.\label{eq:sum_for_touching_0-2}
\end{align}
Analogously  to \eqref{eq:partition_along_0},
\begin{align}
\notag Z_{n_i,n_{i-1}}(0,y)\ge & e^{M_0 r}e^{\sum_{m=n_i+r+1}^{k-1} \xi_m} e^{-M_0}\\ &\times Z_{k,l}(0,0)
e^{\sum_{m=l+1}^{n_{i-1}-r}\xi_m} e^{M_0(r-|y|)}e^{-M_1|y|}.\label{eq:partition_along_0-2}
\end{align}
Dividing \eqref{eq:sum_for_touching_0-2} by \eqref{eq:partition_along_0-2} and taking into account \eqref{eq:lln-consequence-1}, \eqref{eq:lln-consequence-2}, we obtain:
\begin{align}
 \frac{Z_{n_i,n_{i-1}}^{k,l}(y)}{Z_{n_i,n_{i-1}}(0,y)}&\le e^{(3M_1-M_0)r}(2d+1)^{2r}e^{2M_0}e^{-\lambda((k-1-(n_i+r))+(n_{i-1}-r-l))}\notag\\
&\le e^{-2\lambda_0 r}e^{2M_0}e^{-(\lambda(k-1-(n_i+r))+(n_{i-1}-r-l))}\label{eq:generic_frac}.
\end{align}
For $\bar Z_{n_i,n_{i-1}}(y)$ we have
\begin{align*}
\bar Z_{n_i,n_{i-1}}(y)\le &(2d+1)^{r}e^{M_1r}(2d+1)^{n_{i-1}-n_{i}-2r}e^{M_1(n_{i-1}-n_{i}-2r)}(2d+1)^{r}e^{M_1r}
\end{align*}
and
\begin{align*}
Z_{n_i,n_{i-1}}(y)\ge &e^{M_0 r}e^{\sum_{m=n_i+r+1}^{n_{i-1}-r} \xi_m} e^{M_0(r-|y|)}e^{-M_1|y|},
\end{align*}
so that
\begin{equation}
\frac{\bar Z_{n_i,n_{i-1}}(y)}{Z_{n_i,n_{i-1}}(0,y)}\le (2d+1)^{2r}e^{(3M_1-M_0) r} e^{-\lambda(n_{i-1}-n_{i}-2r)}\le e^{-2\lambda_0 r}.
\label{eq:special-frac}
\end{equation}
Now the lemma follows from \eqref{eq:decompositionwrttoching zero2} and \eqref{eq:generic_frac},\eqref{eq:special-frac} 
\epf

\bpf[Proof of Lemma~\ref{lm:terminal}]
Proceeding as above, one can use the fact that all paths contributing to $\hat T_r^{n_i,0}\phi(y)$ do not visit $0$
at first $r$ steps, and 
show that for some $K_9(\lambda)$ and sufficiently large $i$,
\[
\hat T_r^{n_i,0}\phi(y)\le  K_9(\lambda)(2d+1)^re^{(M_1-M_0)r} T_r^{0,n}\phi(y), 
\]
whenever $|y|\le r$. 
Since \eqref{eq:hilbert-metric} implies:
\begin{align*}
&\rho_r(T^{n_i,0}\phi, T^{n_i,0}\phi)\\ 
\le& \ln\left(\max_{|y|\le r} 
\frac{T_r^{n_i,0}\phi(y)\left(1+\frac{\hat T_r^{n_i,0}\phi(y)}{T_r^{n_i,0}\phi(y)}\right)}{T_r^{n_i,0}\psi(y)}\cdot 
\max_{|y|\le r}\frac{T_r^{n_i,0}\psi(y)\left(1+\frac{\hat T_r^{n_i,0}\psi(y)}{T_r^{n_i,0}\psi(y)}\right)}{T_r^{n_i,0}\phi(y)}\right)
\\\le& \rho_r(\phi,\psi)+\max_{|y|\le r} \frac{\hat T_r^{n_i,0}\psi(y)}{T_r^{n_i,0}\psi(y)}\cdot \max_{|y|\le r} \frac{\hat T_r^{n_i,0}\phi(y)}{T_r^{n_i,0}\phi(y)},
\end{align*}
the lemma follows with $K_5(\lambda)=K_9(\lambda)^2$.
\epf


\bpf[Proof of Lemma~\ref{lm:concentration_for_Gibbs_marginals}]
We can write
\[
 \mu_{n_1,n_2}\{|\alpha_{n}|>r\}=\sum_{x_1,x_2\in B_r}\mu_{n_1,n_2}\{\alpha_{n_1}=x_1,\alpha_{n_2}=x_2\}\frac{\tilde Z_{n_1,n_2}(x_1,x_2)}{Z_{n_1,n_2}(x_1,x_2)},
\]
where $\tilde Z_{n_1,n_2}(x_1,x_2)$ denotes the partition function over all paths $\gamma$ defined on $[n_1,n_2]$ with $\gamma(n_1)=x_1$,
$\gamma(n_2)=x_2$, and such that $|\gamma(n)|>r$.
For any $x_1,x_2\in B_r$,
\[
\frac{\tilde Z_{n_1,n_2}(x_1,x_2)}{Z_{n_1,n_2}(x_1,x_2)}\le \sum_{\substack{k_1,k_2:\\|x_1|\le n_1+k_1 \le n-r\\ n+r\le n_1+k_2 \le n_2-|x_2|}}
\frac{\tilde Z^{k_1,k_2}_{n_1,n_2}(x_1,x_2)}{Z^{k_1,k_2}_{n_1,n_2}(x_1,x_2)}+ \frac{\bar Z_{n_1,n_2}(x_1,x_2)}{Z_{n_1,n_2}(x_1,x_2)},
\]
where $\tilde Z^{k_1,k_2}_{n_1,n_2}(x_1,x_2)$ and $Z^{k_1,k_2}_{n_1,n_2}(x_1,x_2)$ are partition functions taken over paths $\gamma$
contributing to $\tilde Z_{n_1,n_2}(x_1,x_2)$ and $Z_{n_1,n_2}(x_1,x_2)$ respectively, with the following restriction:
\begin{align*}
 \sup\{k\le n:\ \gamma_k=0\}=k_1,\\
 \inf\{k\ge n:\ \gamma_k=0\}=k_2,
\end{align*}
and $\bar Z_{n_1,n_2}(x_1,x_2)$ is defined as the partition function over paths contributing to $\tilde Z_{n_1,n_2}(x_1,x_2)$ and never visiting the
origin between $n_1$ and $n_2$. Using Lemma~\ref{lm:main_localization_lemma} to estimate the contribution of  the optimal path $\gamma^*$ 
to the denominator,
we can write:
\begin{align*}
\frac{\tilde Z^{k_1,k_2}_{n_1,n_2}(x_1,x_2)}{Z^{k_1,k_2}_{n_1,n_2}(x_1,x_2)}
&\le \frac{e^{M_1(k_2-k_1)}(2d+1)^{k_2-k_1}}{(2d+1)^{k_2-k_1}
e^{M_1(k_2-k_1)}e^{\lambda(k_2-k_1)}e^{-(M_0+M_1)}}
\\&\le e^{-(M_0+M_1)} e^{-\lambda(k_2-k_1)}. 
\end{align*}
For the last term, we get:
\begin{align*}
 \frac{\bar Z_{n_1,n_2}(x_1,x_2)}{Z_{n_1,n_2}(x_1,x_2)}\le&\frac{e^{M_1(n_2-n_1)}(2d+1)^{n_2-n_1}}{e^{-M_1|x_1|}e^{-(M_0+M_1)}
(2d+1)^{n_2-n_1-|x_1|-|x_2|}}
\\ &\times \frac{1}{e^{M_1(n_2-n_1-|x_1|-|x_2|)}e^{\lambda(n_2-n_1-|x_1|-|x_2|)}e^{-M_1|x_2|}}\\
\le& (2d+1)^{2r} e^{4M_1r}e^{M_0+M_1}e^{-\lambda(n_2-n_1-2r)}.
\end{align*}
and the lemma follows by combining the estimates above.\epf

\section{Acknowledgements} 
The authors would like to thank Leonid Koralov for reading the manuscript and suggesting several useful corrections.
The research of Yuri Bakhtin is partially supported by NSF through CAREER grant DMS-0742424. The research of Konstantin Khanin is partially
supported by NSERC.
\bibliographystyle{alpha}
\bibliography{localiz}

\end{document}